\documentclass{amsart}
\usepackage{graphicx}
\newtheorem{thm}{Theorem}[section]
\newtheorem{cor}{Corollary}[section]
\newtheorem{lem}{Lemma}[section]
\newtheorem{prop}{Proposition}[section]
\theoremstyle{definition}
\newtheorem{defn}{Definition}[section]
\theoremstyle{remark}
\newtheorem{rem}{Remark}[section]
\numberwithin{equation}{section}

\begin{document}

\title{The symmetry axiom in Minkowski planes}%
\author{Jaros\l aw Kosiorek and Andrzej Matra\'{s}}%
\address{Jaros\l aw Kosiorek and Andrzej Matra\'{s}\\
Department of Mathematics and Informatics\\
UWM Olsztyn\\
\.Zo\l nierska 14\\
10-561 Olsztyn\\
Poland}%
\email{matras@uwm.edu.pl}


\dedicatory{Dedicated to H. Karzel}%
\begin{abstract}
The aim of the paper is to give a synthetic proof that in a symmetric Minkowski plane the Benz's (G) axiom holds
(without
using the algebraic representation).\\
\end{abstract}
\maketitle
\section*{Introduction}
The foundations of Minkowski planes were built in the late sixties by R. Artzy, W. Benz, W. Heise, H. Karzel and others
(cf. \cite{A1}, \cite{B1}, \cite{B2}, \cite{HK1}, \cite{HK2}). Two large and well-investigated classes of models of
Minkowski planes are determined by the symmetry axiom (S) (cf. \cite{HK1}) and the rectangle axiom (G) (\cite{B1}).\\
It is well known that the axiom (S) determines the class of miquelian planes over commutative fields and the axiom (G)
the larger class of planes connected with so called Tits near-fields (\cite{Har}). It follows that the symmetry axiom
implies the rectangle axiom but all known proofs have been indirect. They involve the algebraic representation of
the planes. \\
As a generalization of Minkowski planes by omitting the touching axiom (T) hyperbola structures are considered (cf.
\cite{B1}, \cite{HK1}). The geometry of the graphs of a sharply 3-transitive permutation set is a hyperbola structure
(cf. \cite{HK1}). The rectangle axiom is equivalent to the fact that the permutation set describing the hyperbola
structure is closed under composition (cf. \cite{B1}). Hyperbola structures fulfilling (S)
correspond to symmetric sharply 3-transitive permutation sets.\\
In the fine paper \cite{Ka} Karzel proved that every symmetric sharply 3-transitive permutation set is isomorphic to
$\mathrm{PGL}(2,K)$, where $K$ is a commutative field, hence the corresponding hyperbola structure is a symmetric
Minkowski plane and satisfies the rectangle axiom (G).\\
The synthetic proof that the touching axiom (T) is a consequence of the symmetry axiom is given in \cite{HK1}.\\
The aim of this paper is to present a direct, synthetic proof of the implication  $(S) \rightarrow (G) $  and to
explain the geometric connection between the two axioms. Our proof is long, but this point of view provides a natural
and intrinsic characterization of orthogonality in symmetric Minkowski planes.  In some considerations we use
ideas of Karzel from \cite{Ka} and give them a geometric interpretations.\\
In Section 2 we give a synthetic proof that (S) implies both conditions of Dienst (\cite{D})\\
$(\star)$ Every symmetry with respect to any circle is an automorphism.\\
$(\star\star)$ The composition of three symmetries with respect to circles of a bundle is the symmetry with respect to
a circle (of this bundle).\\
In Section 3 we define Minkowski planes of characteristic two independently of the notion of the characteristic of
translation planes.   For other Minkowski planes we introduce harmonic relation of generators. Then we define and
characterize an involutory automorphism with two pointwise fixed generators called harmonic homology.\\
In Section 4 we obtain some configuration theorems which can be regarded as special cases of the rectangle axiom with
orthogonality of suitable circles in assumptions. Especially the configurations from Lemma 4.1 and Corollary 4.4
appeared useful for a description of ordered Minkowski planes in our next prepared paper. In this section we use the
language of permutation set for a brief exposition of results. However all presented statements have clear
geometric interpretations.\\
\textbf{Acknowledgement} The authors wish to thank the reviewer for the many helpful suggestions.
\section{Notations and basic definitions.}
Let $\mathcal{P}$ be a non empty set of \it points\/ \rm and $\Lambda,\Sigma_{1},\Sigma_{2}$ non empty subsets of the
power set of $\mathcal{P}$ with $\Sigma_{1}\cap\Sigma_{2}=\emptyset$; the elements of $\Lambda$ and
$\Sigma_{1}\cup\Sigma_{2}$ will be called \it circles\/ \rm and \it generators\/ \rm respectively.
$\mathcal{M}:=(\mathcal{P},\Lambda,\Sigma_{1}\cup\Sigma_{2})$ is called \it hyperbola structure\/ \rm if the following
conditions are valid:\\
(H1) For each point $p\in\mathcal{P}$ and $i\in\{1,2\}$ there is exactly one generator $G\in\Sigma_{i}$ with
$p\in G$ (notation $[p]_{i}:=G$).\\
(H2) Any two generators $A,B$ with $A\in\Sigma_{1}$, $B\in\Sigma_{2}$ intersect in exactly one point.\\
(H3) Every generator intersects every circle in exactly one point.\\
(H4) For any three distinct points $a,b,c\in\mathcal{P}$ with $[a]_{i}\neq[b]_{i}\neq[c]_{i}\neq[a]_{i}$ for
$i\in\{1,2\}$ there is exactly one circle $C\in\Lambda$ such that $a,b,c\in C$.\\
(H5) There exists a circle containing at least three points.\\
For any point $p\in\mathcal{P}$ the incidence structure $\mathcal{M}^{p}:= (\mathcal{P}^{p},\Lambda^{p})$, where $
\mathcal{P}^{p}:= \mathcal{P}\setminus ([p]_{1}\cup[p]_{2})$ and $ \Lambda^{p}:=
 \{ K \setminus \{ p \} \mid K \in \Lambda, p \in K \} \cup \{ A \setminus [p]_{2} \mid A \in \Sigma_{1}, p \notin A \} \cup \{B \setminus [p]_{1} \mid | B \in \Sigma_{2},p \notin B \}$
is called  \it the derived structure of $\mathcal{M}$ in the point $p$\/ \rm (cf. \cite{HK1}).\\
A hyperbola structure $\mathcal{M}=(\mathcal{P},\Lambda,\Sigma_{1}\cup\Sigma_{2})$ is called a \it Minkowski plane\/
\rm if it has the property:\\
(T) (\it Touching axiom\/ \rm) For any circle $C\in\Lambda$, any $a\in C$ and any $b\in
P\setminus(C\cup[a]_{1}\cup[a]_{2})$ there is exactly one circle $B\in\Lambda$ with $b\in B$ and $C\cap B=\{a\}$.\\
$\mathcal{M}$ is a Minkowski plane if and only if the derived structure $\mathcal{M}^{p}$ in any point $p$ is an affine
plane.\\
For points $p,q$ and $K \in \Lambda $ we write $[p]=[p]_{1} \cup[p]_{2}$, $pq = [p]_{1} \cap [q]_{2}, pK = K
\cap[p]_{1},Kp~=~K \cap[p]_{2}$. For $p,q,r\in\mathcal{P}$ such that $p\notin[q], q\notin[r], p\notin[r]$ the unique
circle through $p,q,r$ is denoted by $(p,q,r)^{\circ}$. A set of points $X$ is called \it concyclic \/ \rm if there
exists a circle $K\in\Lambda$ such that $X\subset K$. If $ p \notin [q]$, the set $ \langle p,q \rangle : = \{ K \in
\Lambda | p,q \in K \}$ is called the $bundle$ with vertices $p,q$. If $ p \in K $, for some circle $K$, the set $(p,K)
= \{ L \in \Lambda | p \in L, L \cap K = \{ p \} \} $ is called the \it pencil\/ \rm determined by the vertex $p$ and
the circle $K$(the pencil of tangent circles with vertex $p$ and direction $K$). Two points $p,q$ are called \it
symmetric with respect to a circle K\/ \rm if $pq$,$qp \in K$. The point symmetric to $p$ with respect to a circle $K$
is denoted by $S_{K}(p)$. For a circle $K$ the involutory bijection $S_{K}$ fixing the points of $K$ and exchanging the
sets $ \Sigma_{1}$, $ \Sigma_{2}$ is called the \it circle symmetry with respect to K.\/ \rm
Two different circles $K,L$ are called \it orthogonal\/ \rm if $S_{K}(p)=p$ for every point $p\in L$(notation: $K\bot L$).\\
In this paper we consider only hyperbola structures satisfying the following symmetry axiom (cf. \cite{B1}, \cite{B2}, \cite{HK1}):\\
(S) For any two circles $K,L$ and a point $p \in L \setminus K$ if $S_{K}(p)\in L$, then $K \bot L$.\\
If a hyperbola structure satisfies (S), then it is a Minkowski plane (cf. \cite{HK1}), so we call it the \it symmetric
Minkowski plane.\/ \rm Until further notice we assume that $\mathcal{M}:=(\mathcal{P},\Lambda,\Sigma_{1}\cup\Sigma_{2})$ is a symmetric Minkowski plane. We repeat this assumption only in the theorems.\\
The main idea of the paper is to present a synthetic proof that every symmetric Minkowski plane fulfils the following rectangle axiom (cf. \cite{B1}, \cite{B2})\\
(G) Let $\{ p_{i} | 1 \leq i \leq 4 \}$, $\{q_{i} | 1 \leq i \leq 4 \}$ be two concyclic sets of different points such
that $\{p_{i}q_{i} | 1 \leq i \leq 4\}$ is also a concyclic set of  points. Then the set of points
$\{q_{i}p_{i}|i=1,...,4\}$ is concyclic.\\
Since the proof of the implication $(\mathrm{S})\rightarrow(\mathrm{T})$ given in \cite{HK1} is synthetic, next we will use the touching axiom.\\

\section{Properties of symmetries with respect to a circle}
We use the following connection between orthogonality and tangency in the pencil of tangent circles with vertex $p$.
\begin{prop}\label{pr2.1}
Let $K,\; L,\; M$ be different circles such that  $p\in K,L,M$, $K \bot M$. Then $L \bot K\Leftrightarrow L \cap M =
\{p\}$.
\end{prop}
\begin{proof}
$\Rightarrow$ See \cite{HK1}, Satz 6.\\
$\Leftarrow$ Let $x\in L$, $x\neq p$ and $x'=S_{K}(x)$. The circle $L'=(x,p,x')^{\circ}$ is orthogonal to $K$ and
tangent to $M$ at $p$ by $(\Rightarrow)$. We obtain $L'=L$ by the touching axiom.
\end{proof}
In order to make all the constructions possible we assume that every circle contains at least 8 points. This is no loss
of generality because in another case the plane is easy to describe (cf. \cite{HS}).
\begin{prop}\label{pr2.2}
If $I,J \in \Lambda$ and $I \cap J = \{ p \}$, then $S_{I}(J)$ is a circle.
\end{prop}
\begin{figure}[h]
\includegraphics[width=0.4\textwidth]{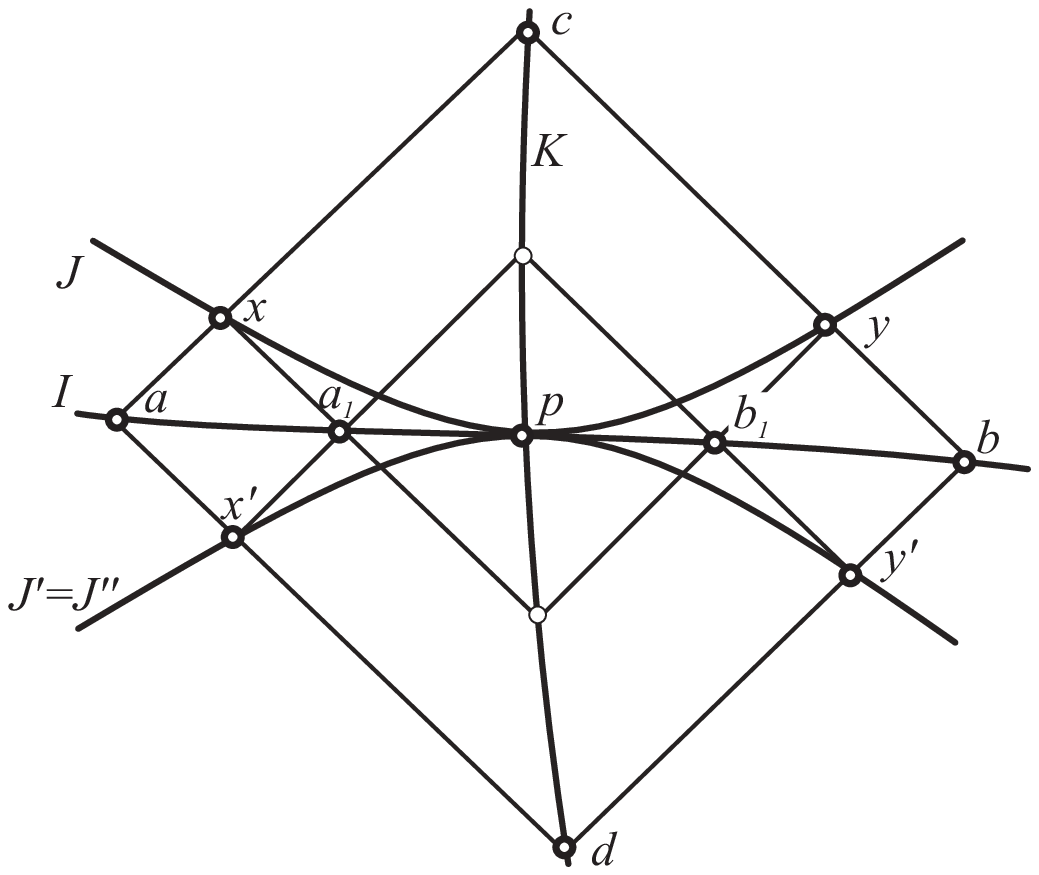}
\caption{}
\end{figure}
\begin{proof}
(fig. 1) Let $x \in J \setminus {p}$, $x' = S_{I}(x)$, and let $J'$ be the unique circle passing through $x'$ and
tangent
to $J$ at $p$. We show that $y' = S_{I}(y) \in J'$ for any $y \in J$.\\
Let $a = xI$, $b = Iy$. Suppose first that $a\neq b$.\\
Let $c = ab$, $d = ba$ and $K = (c,p,d)^{\circ}$. By (S), $K \bot I$. Since $I \cap J = \{ p \}$, we have $K~\bot~J$ by
Proposition~\ref{pr2.1} and the points $x, y$ are symmetric with respect to $K$. If $a_{1}~=~Ix$ and $b_{1}~=~yI$, then
$Ka_{1}~=~Kx~=~yK~=~b_{1}K$ and so $a_{1}, b_{1}$ are symmetric with respect to $K$. Since
$x'K~=~a_{1}K~=~Kb_{1}~=~Ky'$ and $Kx'~=~Ka~=~bK~=~y'K$, it follows that $x', y'$ are symmetric with respect to $K$.
Consider the circle $J''~=~(x',p,y')^{\circ}$. By Proposition~\ref{pr2.1}, $J''~\cap~I~=~\{ p \}$ because $J''~\bot~K$.
According to the touching axiom, $J''~=~J'$ and $y'~\in~J'$.\\
If $a=b$, it is enough to consider the point $z\in J$ such that
$Iz\neq a$ and $zI\neq Iy$. By the first part of the proof we
conclude that $S_{I}(z)\in J'$ and finally that $y'\in J'$.
\end{proof}

\begin{prop}\label{pr2.3}
If $I,J \in \Lambda$ and $I \cap J= \{ a,b \}$ ($a \neq b$), then $S_{I}(J)$ is a circle.
\end{prop}
\begin{proof}
(fig. 2) For any point $x \in J$, $x \neq a,b$, set $x'=S_{I}(x)$. Let $J'=(a,b,x')^{\circ}$ and $y$ be an arbitrary
point such that
$y \in J$, $y \neq a,b,x$. Then $y' = S_{I}(y) \in J'$.\\
Indeed, let $c=ab$, $d=ba$, $e=xy$, $f=yx$ and
$K=(c,d,e)^{\circ}$. By (S), we obtain $K \bot I,J,J'$ and $f \in
K$. We define $e_{1}=eI$, $e_{2}=Ie$, $f_{1}=If$, $f_{2}=fI$. As
$I\bot K$, we deduce that $e_{1}$, $e_{2}$ and $f_{1}$, $f_{2}$
are symmetric with respect to $K$. So we have
$Ky'~=~Kf_{2}~=~f_{1}K~=~x'K$ and $y'K~=~e_{2}K~=~Ke_{1}~=~Kx'$.
Thus $y'$ is symmetric to $x'$ with respect to $K$ and by (S), $y'
\in J'$.
\end{proof}
\begin{figure}[h]
\includegraphics[width=0.5\textwidth]{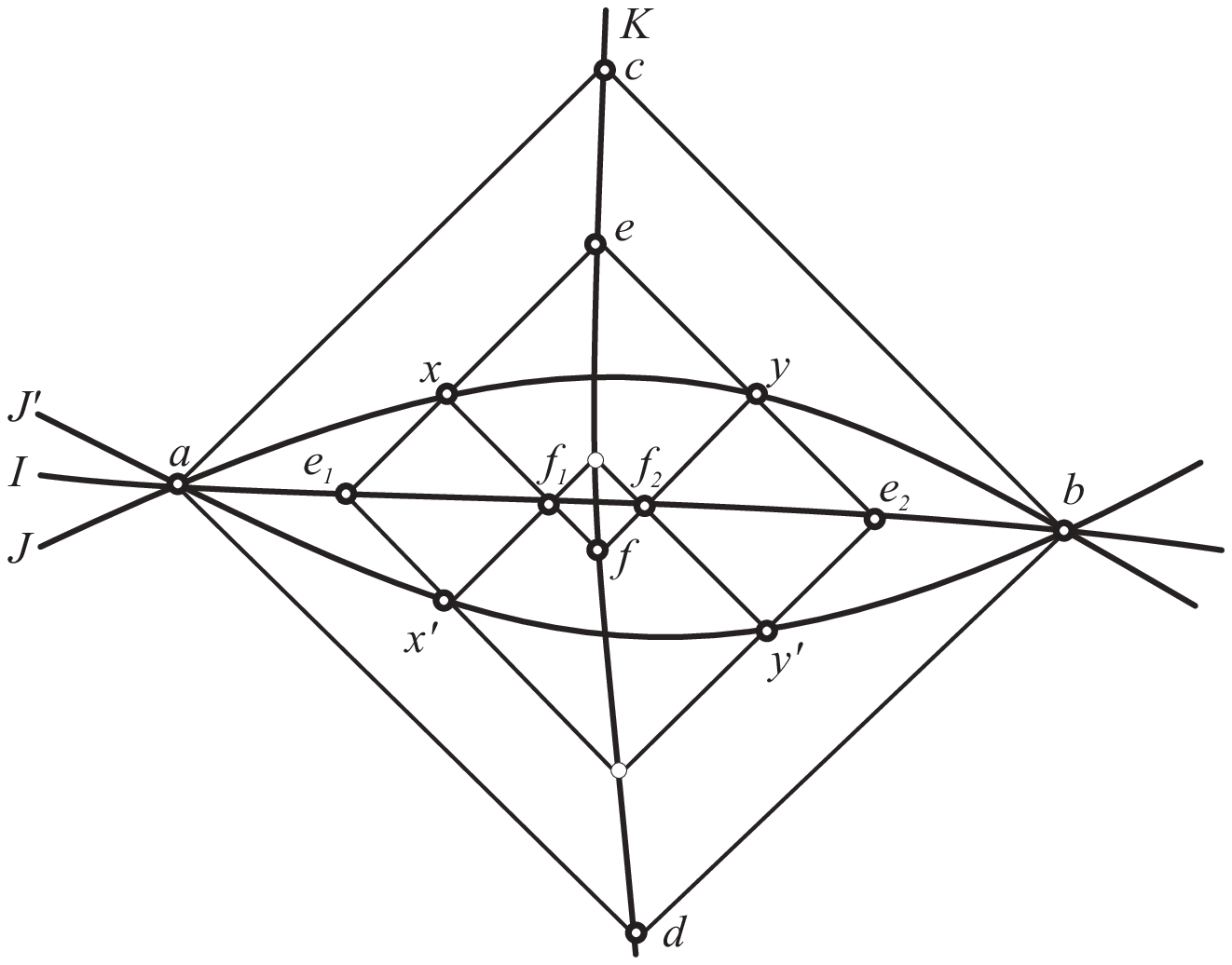}
\caption{}
\end{figure}
\begin{lem} \label{l2.1}
Let $I,J \in \Lambda$ and $x,y \in J \setminus I$ be different points. If $x_{1}=J(xI)$ and $x_{2}=J(x_{1}I)$, then the
points $S_{I}(x)$, $S_{I}(y)$, $S_{I}(x_{1})$, $S_{I}(x_{2})$ are concyclic.
\end{lem}
\begin{proof}
(fig. 3) We put $a=Ix$, $b=yI$, $c=ab$, $d=ba$, $e=S_{J}(d)$ and
$K=(c,d,e)^{\circ}$. By (S), $K\bot I$. We set $b_{1}=Iy$,
$f=b_{1}a$, $g=ab_{1}$, $h=S_{J}(f)$ and $L=(f,g,h)^{\circ}$. By
construction, it follows from (S) that $K,L~\bot~I,J$. We
introduce the following points: $y_{1}=b_{1}J$, $b_{2}=Iy_{1}$,
$y_{2}=b_{2}J$, $b_{3}=Iy_{2}$, $a_{1}=xI$, $a_{2}=x_{1}I$,
$a_{3}=x_{2}I$. Then $y_{1}=S_{L}(x)$. The points $a_{1}$ and
$b_{1}$ are symmetric with respect to $K$ because
$a_{1}K=Kb_{1}=e$ and $I \bot K$. Analogously, we show
successively that the pairs of points $(x_{1},y_{1})$,
$(a_{2},b_{2})$, $(x_{2},y_{2})$, $(a_{3},b_{3})$ are symmetric
with respect to $K$. In the same way we show that $(x,y_{1})$,
$(a_{1},b_{2})$, $(x_{1},y_{2})$, $(a_{2},b_{3})$ are symmetric
with respect to $L$. Let $x'=S_{I}(x)$, $y'=S_{I}(y)$,
$x_{i}'=S_{I}(x_{i})$, $y_{i}'=S_{I}(y_{i})$ ($i=1,2$). Then the
pairs $(x',y')$, $(x_{i}',y_{i}')$ ($i=1,2$) are symmetric with
respect to $K$ and the pairs $(x',y_{1}')$, $(x_{1}',y_{2}')$ are
symmetric with respect to $L$. Indeed, for example for the pair $(x_{1}',y_{2}')$ we have the equalities \\
$x_{1}'L=a_{1}L=Lb_{2}=Ly_{2}'$, \\
$Lx_{1}'=La_{2}=b_{3}L=y_{2}'L$, \\
because the pairs $(a_{1},b_{2})$, $(a_{2},b_{3})$ are symmetric with respect to $L$. From these equalities we see
that $x_{1}'$ and $y_{2}'$ are symmetric with respect to $L$.\\
Consider the circle $M=(x',y',y_{1}')^{\circ}$. Then $M \bot K,L$ because $S_{K}(x')=y'$ and $S_{L}(x')=y_{1}'$. Using
the previous symmetries we get consecutively $x_{1}' \in M$ (because $x_{1}'=S_{K}(y_{1}')$), $y_{2}' \in M$ and
$x_{2}'~\in M$.
\end{proof}
\begin{figure}[h]
\includegraphics[width=0.85\textwidth]{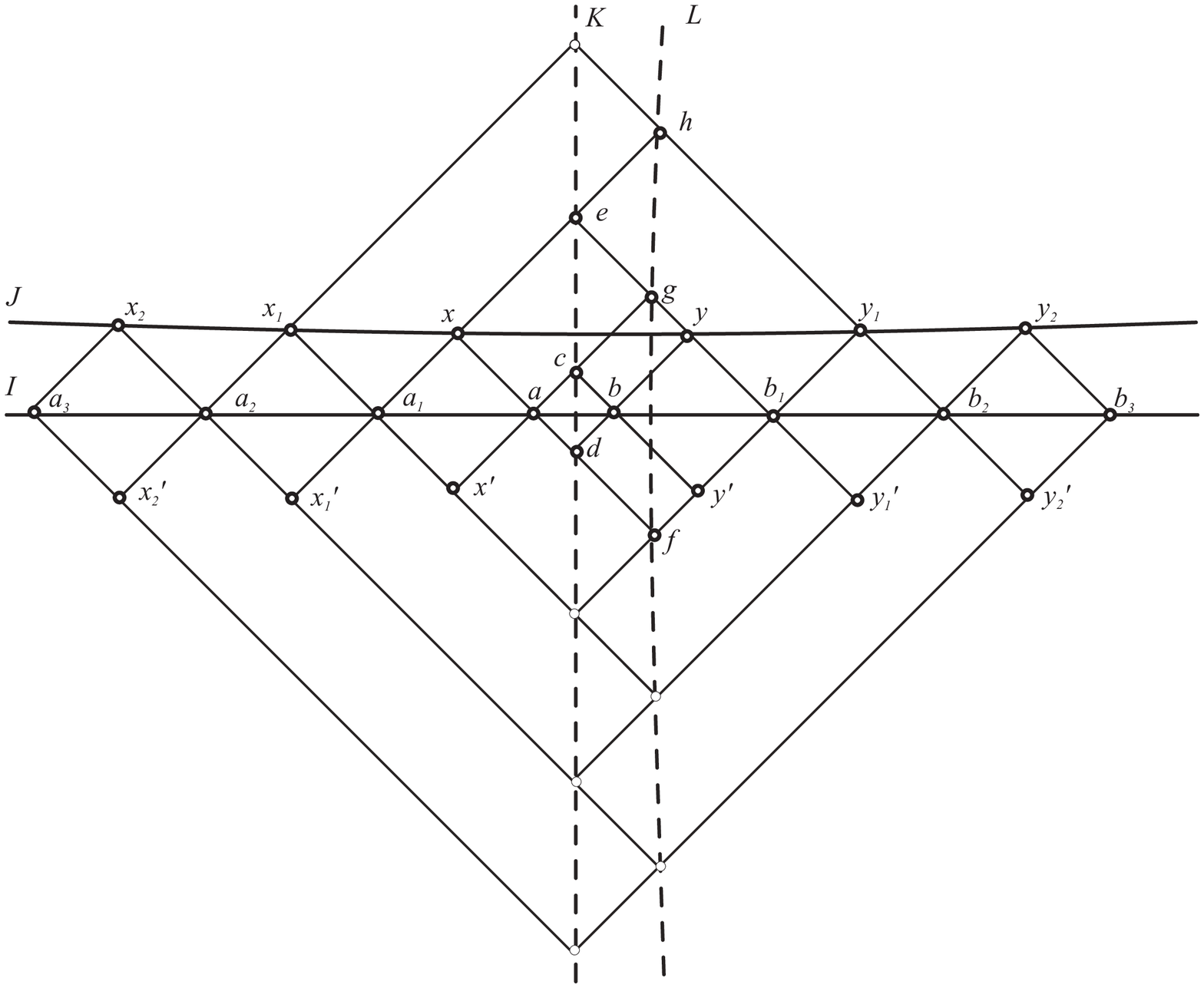}
\caption{}
\end{figure}
\begin{thm} \label{t2.1}
For a symmetric Minkowski plane the symmetry with respect to an arbitrary circle is an automorphism.
\end{thm}
\begin{proof}
Let $J \in \Lambda$. If $I \cap J \neq \emptyset$, then $S_{I}(J)$ is a circle by Propositions~\ref{pr2.2} and
~\ref{pr2.3}. Assume that $I \cap J~=~\emptyset$. Let $x~\in J$ and $x_{1}$, $x_{2}$ be as in Lemma~\ref{l2.1}. Then
Lemma~\ref{l2.1} shows that $S_{I}(y)~\in~(S_{I}(x),~S_{I}(x_{1}),~S_{I}(x_{2}))^{\circ}$ for all $y~\in J$.
\end{proof}
\begin{cor} \label{c2.3}
Every symmetry with respect to a circle preserves orthogonality of circles.
\end{cor}
\begin{thm} \label{t2.2}
Let $\mathcal{M}$ be a symmetric Minkowski plane. If $I,J,K \in \langle p,q\rangle$ for some points $p,~q$, then there
exists $L~\in \langle p,q\rangle$ such that $S_{L}=S_{K}~\circ~S_{J}~\circ~S_{I}$.
\end{thm}
\begin{proof}
(fig. 4) Let $r=pq$, $s=qp$. For every $x \notin [p] \cup [q]$ consider the points $x_{1}=S_{I}(x)$,
$x_{2}=S_{J}(x_{1})$, $x_{3}=S_{K}(x_{2})$ and the circle $M=(r,s,x)^{\circ}$. According to (S), $x_{1},x_{2},x_{3} \in
M$. If $a=xx_{3}$, $b=x_{3}x$, then $a$ and $b$ are symmetric with respect to $M$. For $L=(p,q,a)^{\circ}$, $b \in L$
because $L \bot
M$. We show that for any point $y$, $S_{L}(y)=S_{K} \circ S_{J} \circ S_{I}(y)$.\\
For the points $p, q, r, s, x$ this follows immediately from the construction. Assume that $y \notin [p] \cup [q] \cup
[x]$. Let $y_{1}=S_{I}(y)$, $y_{2}=S_{J}(y_{1})$, $y_{3}=S_{K}(y_{2})$ and $N=(r,s,y)^{\circ}$. Then $y_{1}, y_{2},
y_{3}~\in~N$. If $c=yy_{3}$, $d=y_{3}y$, $z=yx=cb$, then $P=(r,s,z)^{\circ} \bot I, J, K$. It follows that the points
$z_{1}=S_{I}(z)$,
$z_{2}=S_{J}(z_{1})$, $z_{3}=S_{K}(z_{2})$ belong to $P$. We obtain successively the chain of equalities:\\
$z_{1}z=x_{1}x$, $zz_{1}=yy_{1}$, $z_{1}z_{2}=x_{1}x_{2}$, $z_{2}z_{1}=y_{2}y_{1}$, $z_{3}z_{2}=x_{3}x_{2}$,
$z_{2}z_{3}=y_{2}y_{3}$. From the first and fifth equality we get $z_{3}z=b$, and from the second and sixth $zz_{3}=c$,
so $S_{P}(b)=c$.
Analogously, $S_{P}(a)=d$. Since $P \bot L$, it follows that $c,d~\in~L$ and $y_{3}=S_{L}(y)$.\\
If $y\in[x]$, the above arguments apply with some simplifications ($z=x$ or $z=y$).\\
If $y~\in~[p]~\cup~[q]$, then without loss of generality we assume that $y~\in~[p]_{1}, y~\neq~p,r$. For any $a \in
[y]_{2} \setminus [p]$ by the previous part of the proof we get $S_{K}S_{J}S_{I}(a)=S_{L}(a)$. It follows that
$S_{K}S_{J}S_{I}([y]_{2})=S_{L}([y]_{2})$ because $[a]_{2}=[y]_{2}$. We also have
$S_{K}S_{J}S_{I}([y]_{1})=S_{L}([y]_{1})$ because $[y]_{1}=[p]_{1}$, so $S_{K}S_{J}S_{I}(y)=S_{L}(y)$.
\end{proof}
\begin{figure}[h]
\includegraphics[width=0.75\textwidth]{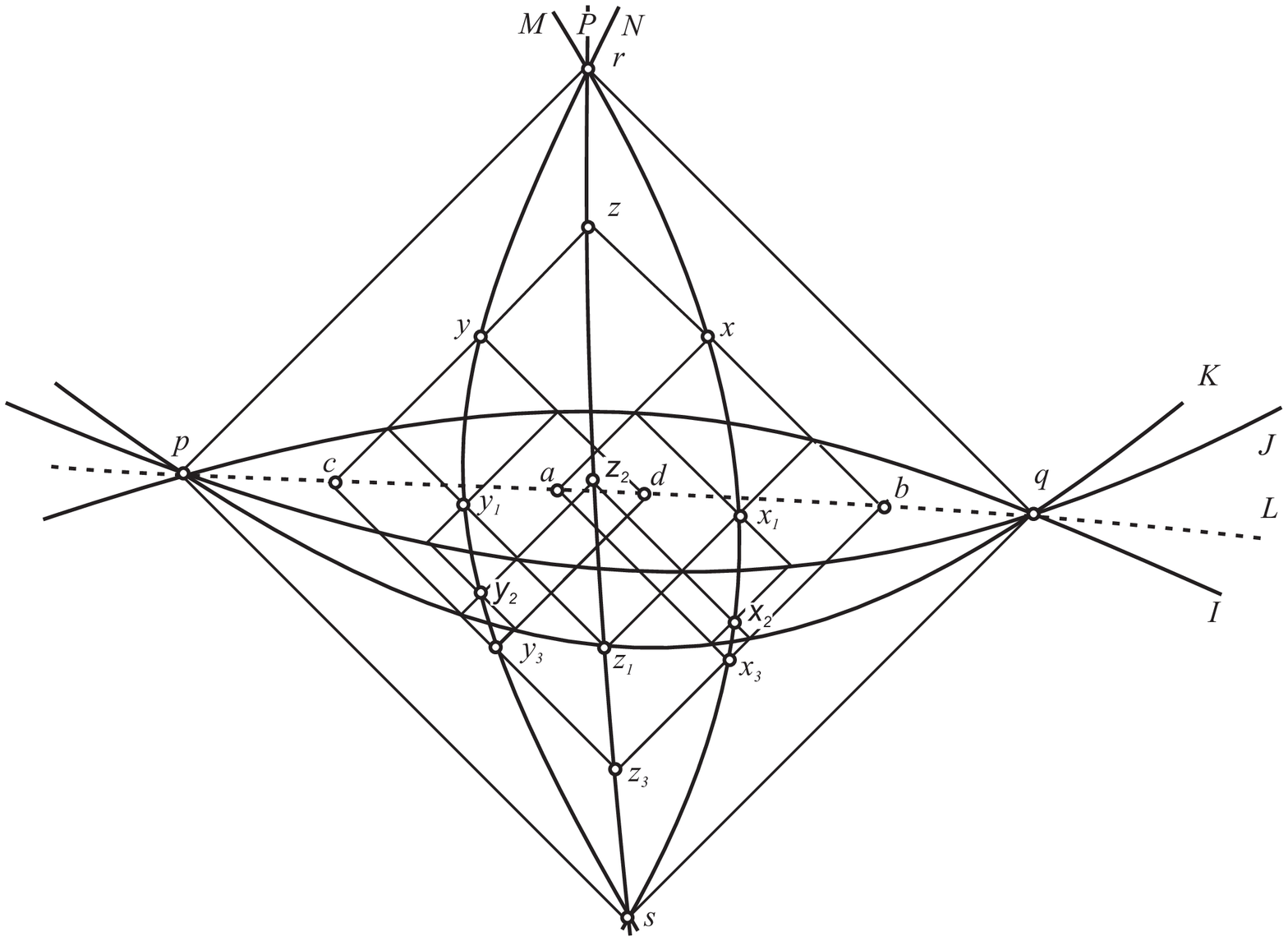}
\caption{}
\end{figure}

Just as for bundles, we have the theorem on reduction for pencils of tangent circles.
\begin{prop}\label{pr2.4}
If $I,J~\in~(p,K)$, then there exists $L~\in~(p,K)$ such that $S_{L}=S_{K}~\circ~S_{J}~\circ~S_{I}$.
\end{prop}
\begin{proof}
The composition $S_{K} \circ S_{J} \circ S_{I}$ is an automorphism which induces on the derived affine plane $
\mathcal{M}^{p}$ a central collineation with improper center corresponding to the ideal point of the pencil of circles
orthogonal to $K$ in $p$. The axis of the collineation is a proper line (because each symmetry exchanges the sets
$\Sigma_{1}$ and $\Sigma_{2}$), so we get a pointwise fixed circle passing through $p$, tangent to $K$.
\end{proof}
\section{Characteristic of a Minkowski plane and harmonic relation of generators of the same kind}
\begin{prop}\label{pr3.1}
If $I,J~\in~\Lambda$, then\\
$I~\bot~J~\Leftrightarrow~S_{J}~\circ~S_{I}$ is an involution.
\end{prop}

\begin{proof}
(fig. 5) $\Rightarrow$ Since $I~\bot~J$, for $x~\in~I~\cup~J$, we have $S_{J}\circ S_{I}(x)=S_{I}\circ S_{J}(x)$. If
$x~\notin~I~\cup~J$,
we set $x'=S_{I}(x)$, $x''=S_{J}(x')$ and $x_{1}=S_{I}(x'')$. We show that $x=S_{J}(x_{1})$.\\
Let $a=xI$, $b=Ix$, $c=x''I$, $d=Ix''$. We obtain $Ja=Jx'=x''J=cJ$, so the points $c$ and $a$ are symmetric with
respect to $J$. Similarly, $b$ and $d$ are symmetric with respect to $J$. We obtain $xJ=aJ=Jc=Jx_{1}$ and $Jx=Jb=dJ=x_{1}J$. Hence $x,x_{1}$ are symmetric with respect to $J$.\\
$\Leftarrow$ For any $x\in I$ we have $S_{I}\circ S_{J}(x)=S_{J}\circ S_{I}(x)=S_{J}(x)$, hence $S_{J}(x)\in I$.
\end{proof}
\begin{figure}[h]
\includegraphics[width=0.3\textwidth]{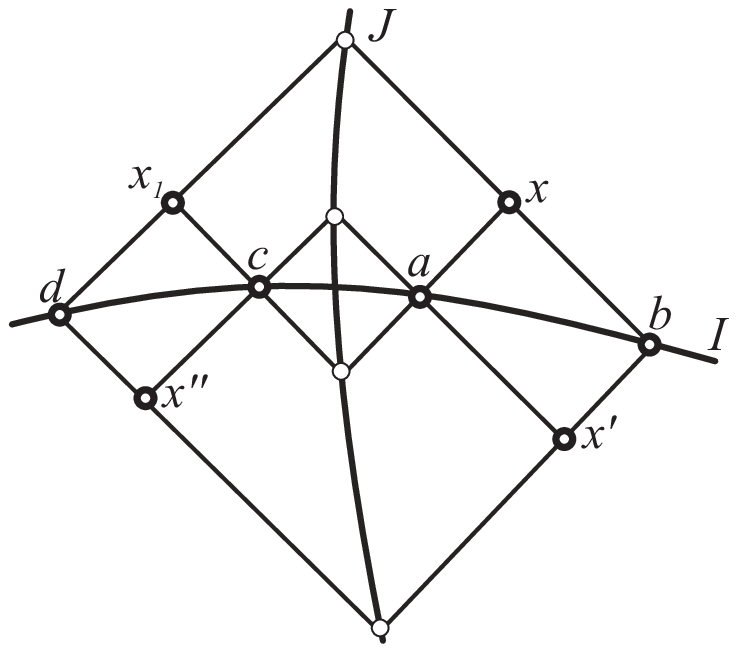}
\caption{}
\end{figure}
\begin{defn}\label{d3.1}
(\cite{J}, Def. 2.2) An automorphism $\varphi$ of a Minkowski plane is called a \it double homothety\/ \rm if there
exist points $p,q$, $(p\notin[q])$ such that $p,q,pq,qp\in \rm fix \varphi$ and $\varphi_{|\mathcal{M}^{p}}$,
$\varphi_{|\mathcal{M}^{pq}}$ are homotheties of affine planes.
\end{defn}
\begin{cor} \label{c3.1}
\rm (\cite{J}, Prop. 2.2, p. 422) \it If $I,J\in\Lambda$, $I\bot$$J$, $I\cap J=\{a,b\}$, then $S_{J}\circ S_{I}$ is a
double homothety.
\end{cor}

\begin{prop}\label{pr3.2}
If $I,J,K\in\Lambda$, $I\bot J,K$, $I\cap J=\{p\}$ and $p\in K$, then $I\cap K=\{p\}$.
\end{prop}
\begin{proof}
Since $I\cap J=\{p\}$ and $I\bot K$, we have $K\bot J$ by Proposition \ref{pr2.1}. As $K,I\bot J$ and $p\in I,J,K$, it
follows that $K\cap I=\{p\}$ by Proposition \ref{pr2.1}.
\end{proof}
\begin{prop}\label{pr3.3}
If $I\cap J=\{p,q\}$ and $I\bot J$, then for every circle $K\in \langle p,q\rangle$ there exists exactly one $L\in
\langle p,q\rangle$ such that $L\bot K$.
\end{prop}
\begin{proof}
By Theorem \ref{t2.2}, there exists $L\in <p,q>$ such that $S_{L}=S_{K}~\circ~S_{J}~\circ~S_{I}$. So we get
$S_{J}~\circ~S_{I}=S_{K}~\circ S_{L}$ and $K \bot L$ by Proposition \ref{pr3.1}. If $L,L' \in <p,q>$, $L,L'\bot K$ and
$L'\neq L$, we obtain a contradiction with Proposition \ref{pr2.1}.
\end{proof}
The following theorem allows us to introduce the definition  of symmetric Minkowski planes of characteristic 2.
\begin{thm} \label{t3.1}
Let $\mathcal{M}$ be a symmetric Minkowski plane. If there exists a pair of orthogonal and tangent circles, then every
pair of tangent circles is orthogonal and every pair of orthogonal circles has at most one common point.
\end{thm}
\begin{proof}
First we prove the second part of the theorem. Since any pair of different points can be transformed in any other by a
 composition of two symmetries and each symmetry preserves orthogonality (Corollary \ref{c2.3}), it is enough to show that the following is
impossible: $I\bot
J$, $I\cap J=\{p\}$ and $K\bot L$, $K,L\in \langle p,q\rangle$ for some $I,J,K,L\in\Lambda$, $p,q\in \mathcal{P}$, $p\notin [q], q\in I$.\\
Let $K'$ be a circle such that $q\in K', K'\cap J=\{p\}$. By Proposition \ref{pr3.3}, there exists a circle $L'\in
\langle p,q\rangle$ with $L'\bot K'$. Then $L'\bot I,J$ and by Proposition \ref{pr2.1}, $L'$ is tangent to $K'$,
contrary to $L'\in <p,q>$.\\
To prove the first part assume that $I\cap J=\{p\}$. If $K\in\Lambda$, $p\in K$, $K\neq I$ and $K\bot J$, then by the
already proven part of the theorem we get $K\cap J=\{p\}$. In particular $K\cap I=\{p\}$ and $K\bot J$. From
Proposition \ref{pr2.1} we obtain $I\bot J$.
\end{proof}
By Theorem \ref{t3.1}, we can introduce the following definition.
\begin{defn}\label{d3.2}
A symmetric Minkowski plane $\mathcal{M}$ is of characteristic 2 (notation $\rm char\mathcal{M}=2$) if there exists a
pair of orthogonal tangent circles.
\end{defn}
We use the notation $\rm char\mathcal{M}\neq 2$ in the case $\mathcal{M}$ is not of characteristic 2 although we don't
define another  characteristic of Minkowski planes. This doesn't arise confusion and agrees with the usually used
notation (cf. \cite{D}, \cite{Har}, \cite{HK1}).\\
From the definitions we remark that $\rm char\mathcal{M}\neq 2$ exactly when for some points $p, q$ with $p\notin [q]$
there exists a double homothety with centers $p, q$ (cf. \cite{J}). For such planes
the existence of one homothety implies that all the possible double homotheties exist.\\
From Proposition \ref{pr3.1} and \cite[Proposition 2.2]{J}, we get the following characterization of involutory
homotheties of symmetric Minkowski planes.
\begin{prop}\label{pr3.4}
If $\varphi$ is an automorphism of a Minkowski plane $\mathcal{M}$ with $\rm char\mathcal{M}\neq 2$, then the following are equivalent:\\
(i) $\varphi$ is an involutory homothety,\\
(ii) $\varphi$ is a double homothety,\\
(iii) $\varphi$ is a superposition of two symmetries with respect of two orthogonal intersecting circles.
\end{prop}
\begin{rem}\label{r3.1}
The equivalence of (i) and (ii) is proved in \cite{JKM} without a further assumption.
\end{rem}
By Proposition \ref{pr3.4} and \cite[Corollary 2]{JKM} if $\rm char\mathcal{M}\neq 2$ for any $a,b$ ($a\notin[b]$)
there exists exactly one double homothety with centers $a,~b$ which will be denoted by $\mathbf{H}_{a,b}$.
\begin{lem}\label{l3.1}
If $\rm char\mathcal{M}\neq 2$ and points $a,b,c,d$ satisfy the conditions $[a]_{i}=[c]_{i}$, $[b]_{i}=[d]_{i}$ for a
fixed $i$ $(i=1,2)$, then $\mathbf{H}_{a,b}(X)=\mathbf{H}_{c,d}(X)$ for every $X\in\Sigma_{i}$ .
\end{lem}
\begin{proof}
First we assume $a=c$. The automorphism
$\varphi=\mathbf{H}_{a,d}\circ\mathbf{H}_{a,b}$, being a
superposition of two homotheties, induces a collineation with
improper axis of the derived affine plane $\mathcal{M}^{a}$ . If
$\varphi'=\varphi_{|\mathcal{M}^{a}}$ had a proper center, then
there would exist $e\in[a]_{i}$, $f\in[b]_{i}\cap[e]_{3-i}$ such
that $\varphi(e)=e$, $\varphi(f)=f$. Then $\varphi$ induces a
homothety with center $f$ on the plane $\mathcal{M}^{a}$.
Analogously $\varphi$
induces a homothety with center $e$ of the plane $\mathcal{M}^{a'}$, where $\{a'\}=[a]_{3-i}\cap[b]_{i}$.\\
Hence, $\varphi=\mathbf{H}_{a,f}$ would be an involution and
$\mathbf{H}_{a,d}(b)$ a fixed point of $\mathbf{H}_{a,b}$, a
contradiction. It follows that the center of $\varphi'$ is an
improper point corresponding to the class $\Sigma_{1}$ and
$\mathbf{H}_{a,b}(X)=\mathbf{H}_{a,d}(X)$ for every $X\in
\Sigma_{i}$ . If $a\neq c$, then by the above, for every
$X\in\Sigma_{i}$ we get
$\mathbf{H}_{a,b}(X)=\mathbf{H}_{c,b}(X)=\mathbf{H}_{c,d}(X)$.
\end{proof}
For $\mathcal{M}=(\mathcal{P},\Lambda,\Sigma_{1}\cup\Sigma_{2})$ with $\mathrm{char}\mathcal{M}\neq 2$,
$A,B\in\Sigma_{i}$ and $A\neq B$ we may put the following two definitions.
\begin{defn}\label{d3.3}
$X,Y\in\Sigma_{i}$ are \it harmonic conjugate\/ \rm with respect to $A,B$ if there exist $a\in A$, $b\in B$ such that
$\mathbf{H}_{a,b}(X)=Y$ and $X,Y\neq A,B$.
\end{defn}

\begin{defn}\label{d3.4}
The \it harmonic homology\/ \rm with axis $A,B$ (denoted by $\mathbf{S}_{A,B}$) is a permutation of
$\mathcal{P}$ defined as follows:\\
(i) for $x\in A\cup B$, $\mathbf{S}_{A,B}(x) =x$\\
(ii) for $x\notin A\cup B$, $\mathbf{S}_{A,B}(x) =y\Leftrightarrow [x]_{3-i}=[y]_{3-i}$ and the generators $[x]_{i}$,
$[y]_{i}$ are harmonic conjugate with respect to $A,B$.
\end{defn}
By Lemma \ref{l3.1} Definition \ref{d3.3} is independent of the points $a,b$.
\begin{cor}\label{c3.2}
For any $A,B\in\Sigma_{i}$ such that $A\neq B$, $\mathbf{S}_{A,B}$ is an involutory automorphism of the Minkowski
plane.
\end{cor}
\begin{proof}
$\mathbf{S}_{A,B}$  is an involutory bijection preserving the sets $\Sigma_{i}$ because double homothety is an
involution. Let $K\in\Lambda$, $a=A\cap K$, $b=B\cap K$. We show that $\mathbf{S}_{A,B}(K)=L$, where $L\in <a,b>$ and
$L\bot K$. By Proposition \ref{pr3.4}, $\mathbf{H}_{a,b}=S_{L}\circ S_{K}$. If $x\in K$, $x'\in L$ and
$[x]_{3-i}=[x']_{3-i}$, then $\mathbf{S}_{A,B}([x]_{i})=S_{L}\circ
S_{K}([x]_{i})=S_{L}([x]_{3-i})=S_{L}([x']_{3-i})=[x']_{i}$ and $x'=\mathbf{S}_{A,B}(x)$.
\end{proof}
\begin{cor}\label{c3.3}
If $\rm char\mathcal{M}\neq 2$, $X,Y\in\Sigma_{i}$, $X\neq Y$, $K,L\in\Lambda$, then\\
$\mathbf{S}_{X,Y}(K)=L\Leftrightarrow L\bot K \wedge K\cap L=K\cap(X\cup Y)$ .
\end{cor}
We now turn to the case of characteristic two.
\begin{prop}\label{p3.5}
Let $\rm char\mathcal{M}=2$. If $X\neq Y$, $X,Y\in\Sigma_{i}$ and $p\notin X\cup Y$, then there exists an involutory
translation exchanging $X$, $Y$ and fixing pointwiese $[p]_{i}$.
\end{prop}
\begin{proof}
Let $K\in \Lambda$ with $p\in K$ and $x,x'$ be different points of $K$ with $x,x'\neq p$. The circle
$L=(p,xx',x'x)^{\circ}$ is orthogonal to $K$. By Definition \ref{d3.2} and Theorem \ref{t3.1} $L$ is tangent to $K$ in
a point $p$. The automorphism $S_{L}\circ S_{K}$, according to Proposition \ref{pr3.1}, is an involution exchanging
$x,x'$, and induces a translation of $\mathcal{M}^{p}$. In the same manner we can obtain the transitivity of
translations of $\mathcal{M}^{p}$ in any other direction determined by a circle.\\
By  \cite[Theorem 4.19, p. 100]{Hu}, we obtain the transitivity of translations with the direction determined by
$\Sigma_{3-i}$. From the proof of \cite[Theorem 4.19]{Hu} it follows that these translations are compositions of
translations with directions determined by circles, so they are induced by automorphisms of $\mathcal{M}$. They are
also involutions by \cite[Theorem 4.14, p. 97]{Hu}.
\end{proof}
\begin{rem}\label{r3.1}
From Corollary \ref{c3.2} and Proposition \ref{p3.5} it follows that on a symmetric Minkowski plane for different
generators $A,X,Y\in\Sigma_{i}$ there exists an involutory automorphism which fixes $A$ pointwise  and exchanges $X$
and $Y$ (see \cite{Har} for this fact for Minkowski planes with axiom (G)).
\end{rem}
\begin{prop}\label{pr3.6}
Let $\rm char\mathcal{M}=2$. If circles $K,L$ are tangent at $p$ and $M\bot K,L$, then $M\in(p,K)$
\end{prop}
\begin{proof}
If $p\notin M$, we obtain a contradiction with Theorem \ref{t3.1} because the points $p,S_{M}(p)$ are two different
common points of the orthogonal circles $K, L$.
\end{proof}
\section{The proof of postulate G}
We begin with the main lemma.
\begin{lem}\label{l4.1}
Let $p,q\in \mathcal{P}$ and $K,L,M,N\in\Lambda$  satisfy the following conditions:\\
(i) $p\neq q$, $p\in [q]_{2}$, $p\in K\cap L$, $q\in M\cap N$,\\
(ii) $K,M\bot L,N$.\\
Then the conditions $k\in K$, $l\in L$, $m\in M$, $n\in N$ and $k\in[l]_{1}$, $l\in[m]_{2}$, $m\in[n]_{1}$ imply
$k\in[n]_{2}$.
\end{lem}
\begin{proof}
(fig. 6) Assume that $char\mathcal{M}\neq 2$. We define the points $r=pM$, $s=pN$, $t=qL$ and $u=qK$. By
orthogonalities of suitable circles, we get $r\in[t]_{2}$ and $s\in[u]_{2}$. For $x\in K\cap L$, $x\neq p$ consider the
harmonic homology $\varphi=\mathbf{S}_{P,X}$, where $P=[p]_{2}$, $X=[x]_{2}$. By definition of $\varphi$, we get
$\varphi(K)=L$ and $\varphi(t)=u$. Hence $\varphi(s)=r$ and $\varphi(M)=N$. By assumptions, $\varphi(l)=k$,
$\varphi(m)=n$ and
$l\in[m]_{2}$. We finally  get $k\in [n]_{2}$.\\
In the case $char\mathcal{M}=2$ the proof is analogous. Instead of the automorphism $\mathbf{S}_{P,X}$ we use a
translation with pointwise fixed generator $P$ which maps $t$ to $u$; such a translation exists by Proposition
\ref{p3.5} and is involutory.
\end{proof}
\begin{figure}[h]
\includegraphics[width=0.4\textwidth]{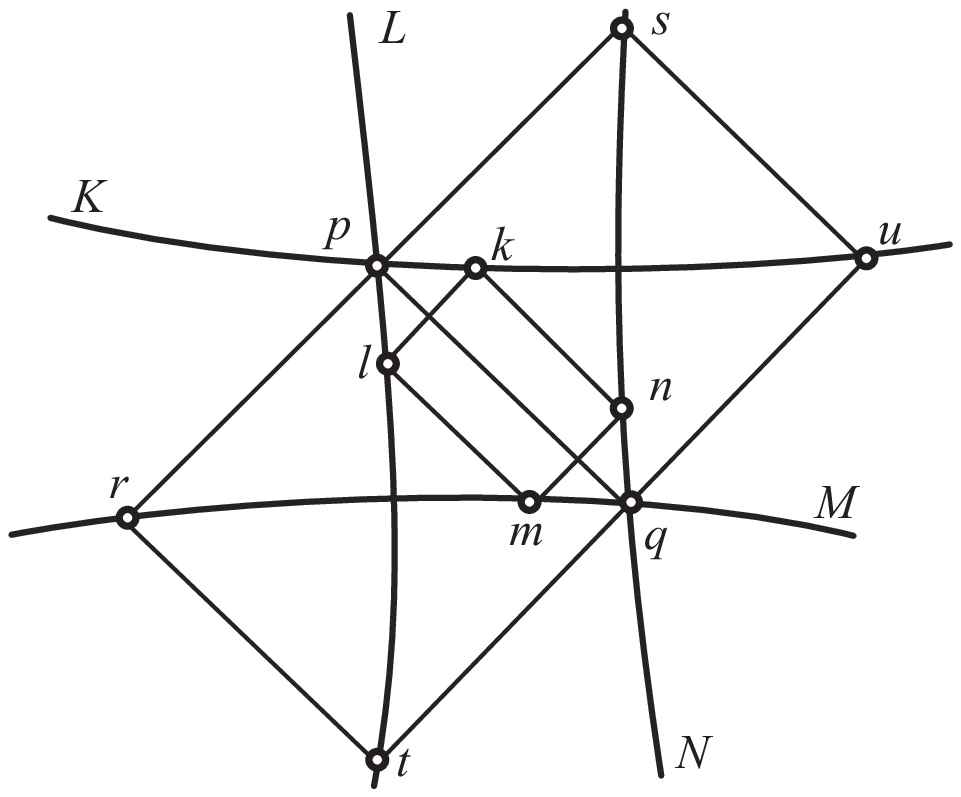}
\caption{}
\end{figure}
\begin{lem}\label{l4.2}
Let $p,r,s,q,a,b,c,d$ be different points and $K,L,M,N$ be different circles such that $s=pr, q=rp, b=ca, d=ac$,
$K,M\in\langle p,r\rangle$,
$L,N\in\langle q,s\rangle$, $a\in L$, $b\in M$, $c\in N$, $d\in K$.\\
If $a'\in L$, $b'\in M$, $c'\in N$ and $b'=c'a'$, $d'=a'c'$, then $d'\in K$.
\end{lem}
\begin{proof}
(fig. 7) Consider the points $t=dd'$, $u=bb'$ and the orthogonal circles $T=(s,t,q)^{\circ}$, $U=(u,p,r)^{\circ}$. If
$u'=\mathbf{S}_{N}(u)$, $u''=\mathbf{S}_{L}(u)$, then $u',u''\in U$ because  $L,N\bot U$. Hence
$\mathbf{S}_{T}(u')=u''$ and $\mathbf{S}_{T}(d)=d'$. It follows that $d'\in K$ because $T\bot K$.
\end{proof}
\begin{figure}[h]
\includegraphics[width=0.6\textwidth]{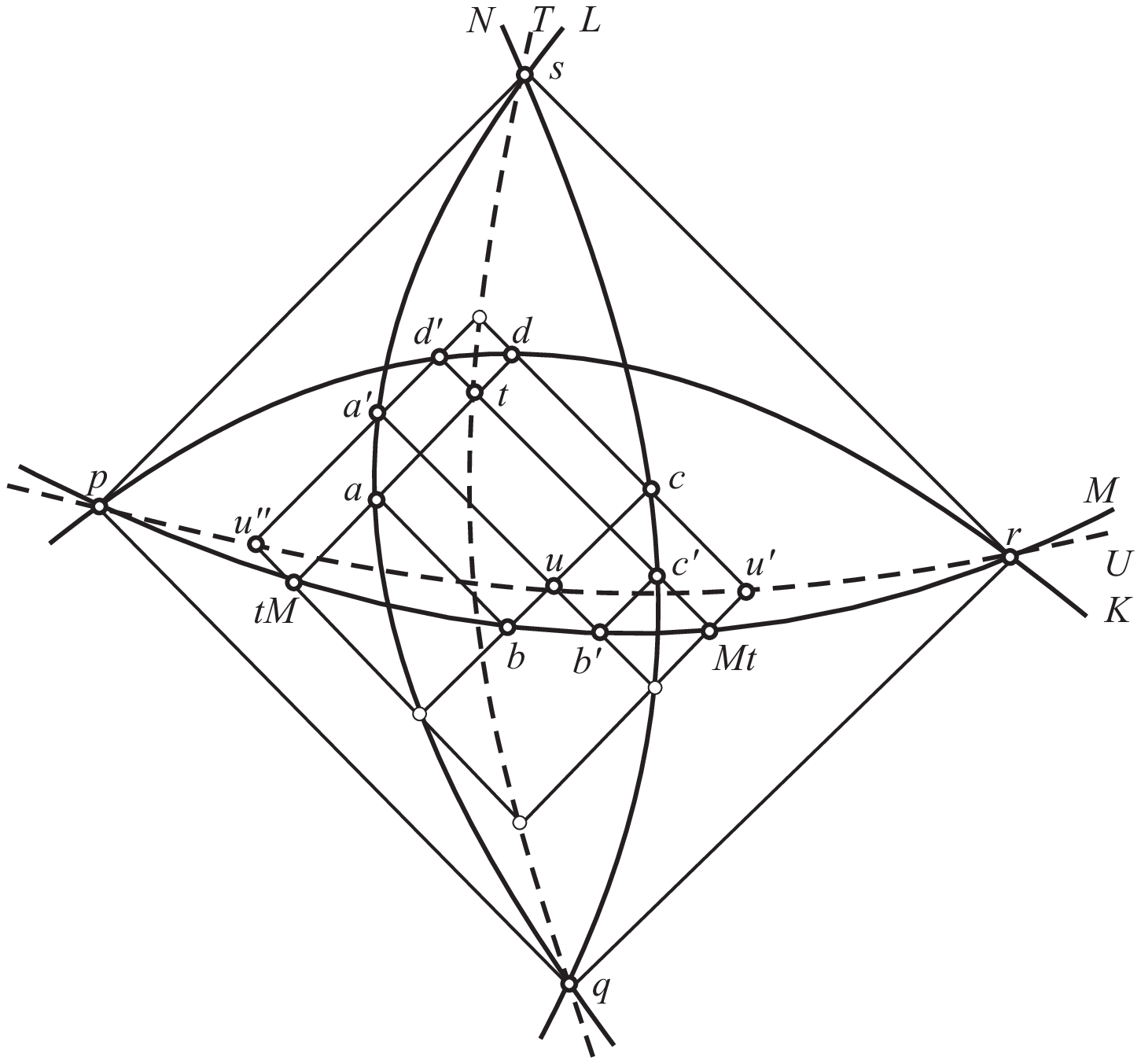}
\caption{}
\end{figure}
In the remainder of this section we fix an arbitrary circle $E$.\\
It is well known that every $K\in\Lambda$ determines the permutation $f_{K}: E \rightarrow E$ by the formula
$f_{K}(x)=(Kx)E$ for $x\in E$ and vice versa (cf. \cite{B1}, \cite{Har}, \cite{HK1}, \cite{Ka}). By definition,
$f_{E}=id_{E}$. According to (S), if $K\bot E$, then $f_{K}^{2}=id_{E}$. The involution corresponding to a circle $K$,
$K\bot E$ will be denoted by $i_{K}$. According to Theorem \ref{t2.1}, for any $f_{K}$, there exists the permutation
$f_{K}^{-1}:=f_{K'}$, where $K'=S_{E}(K)$.\\
Now we are going to show that the permutation set $\{f_{K}|K\in\Lambda\}$ is closed under composition. This is
equivalent to the rectangle axiom.
\begin{cor}\label{c4.1}
$\mathrm{(i)}$ If $L,N\in\Lambda,~q\in E\cap(N\setminus L)$ and $L,N\bot E$, then there exists $K\in\Lambda$ such that
$Lq\in K$,
$K\bot L,N$ and $f_{K}=i_{L}\circ i_{N}.$\\
$\mathrm{(ii)}$ If $\mathrm{char}\mathcal{M}\neq 2$ or $K\bot\!\!\!\!\angle E$, then for every $q\in E\setminus K$
there exist circles $L,N$ such that $f_{K}=i_{L}\circ i_{N}$, $q\in N$, $Kq\in L$ and $L,N\bot K$.
\end{cor}
\begin{proof}
Assuming $E:=M$ in Lemma \ref{l4.1}, we get $f_{K}=i_{L}\circ i_{N}$, so it is enough to show the existence of suitable circles.\\
(i) If $L\bot\!\!\!\!\angle N$, then we set $K:=(qL,S_{N}(qL),S_{L}(S_{N}(qL)))^{\circ}$. If $L\bot N$ and
$\mathrm{char}\mathcal{M}\neq 2$ the unique circle $K$ such that $K\in \langle qL,S_{N}(qL) \rangle$ and $K\bot L$
(Proposition \ref{pr3.3}) is the circle we are looking for. The case $\mathrm{char}\mathcal{M}=2$ and $L\bot N$
is impossible by Proposition \ref{pr3.6}.\\
(ii) In the case $K\bot\!\!\!\!\angle E$ the assertion follows for the circles $N=(q,S_{K}(q),S_{E}(S_{K}(q)))^{\circ}$
and $L=(qK,S_{E}(qK),S_{K}(S_{E}(qK)))^{\circ}$. If $K\bot E$ and $\mathrm{char}\mathcal{M}\neq 2$,  we take the
circles $N,L$ such that $N\in \langle q,S_{K}(q) \rangle$,  $N\bot E$, $L\in \langle qK,S_{E}(qK) \rangle$ and $L\bot
K$.
\end{proof}
\begin{rem}\label{r4.1}
Corollary \ref{c4.1} enables one to represent any permutation $f_{K}$ as a composition of two involutions one of which
has a fixed point (in the case $char\mathcal{M}= 2$ with the assumption that $f_{K}$ is not an involution).
\end{rem}
\begin{rem}\label{r4.2}
In Lemma \ref{l4.1} $K\bot M$ (hence $N\bot L)$ is possible only in the case $char\mathcal{M}\neq2$.
\end{rem}
Analogously to Corollary \ref{c4.1}, from Lemma \ref{l4.2}, we obtain
\begin{cor}\label{c4.1a}
$\mathrm{(i)}$ If $L,N\in\langle q,s\rangle$, $L,N\bot E$ and $E\notin \langle q,s \rangle$, then there exists
$K\in\Lambda$ such
that $K\in \langle qE,Eq \rangle$ and $f_{K}=i_{L}\circ i_{N}$.\\
$\mathrm{(ii})$ If $K,E\in\langle p,r \rangle$ and $c\notin[p]\cup[r]$, then there exist circles $L,N$ such that $c\in
N$, $L,N\in \langle pr,rp \rangle$ and $f_{K}=i_{L}\circ i_{N}$.
\end{cor}
\begin{prop}\label{pr4.1}
If $L,N \bot E$ and $p\in E\cap L\cap N$, then there exists a circle $K$ such that $K\cap E=\{p\}$ and $i_{N}\circ
i_{L}=f_{K}$.
\end{prop}
\begin{proof}
Let $char\mathcal{M}\neq2$. There exists a point $q\neq p$ such that $q\in E\cap L$ and $q\notin N$ by Proposition
\ref{pr2.1} (fig. 8). From  part (i) of Corollary \ref{c4.1} it follows that there exists $K\in\Lambda$, $K\bot L,N$
such that $i_{L}\circ i_{N}=f_{K}$. Since $p\in L\cap N$, we have $f_{K}(p)=p$ and $K\cap E=\{p\}$.\\
If $char\mathcal{M}=2$, by Lemma \ref{l2.1} there exists $K\in\Lambda$ such that $S_{K}=S_{L}\circ S_{E}\circ S_{N}$.
If $x\in E$ and $y=Nx$, then $S_{K}(y)=S_{L}S_{E}(y)$ (fig. 9). Hence $(L(yE))K=Ky$ and this means that $i_{L}\circ
i_{N}(x)=f_{K}(x)$. We remark that in this case $f_{K}$ is an involution.
\end{proof}
\begin{figure}[h]
\includegraphics[width=0.7\textwidth]{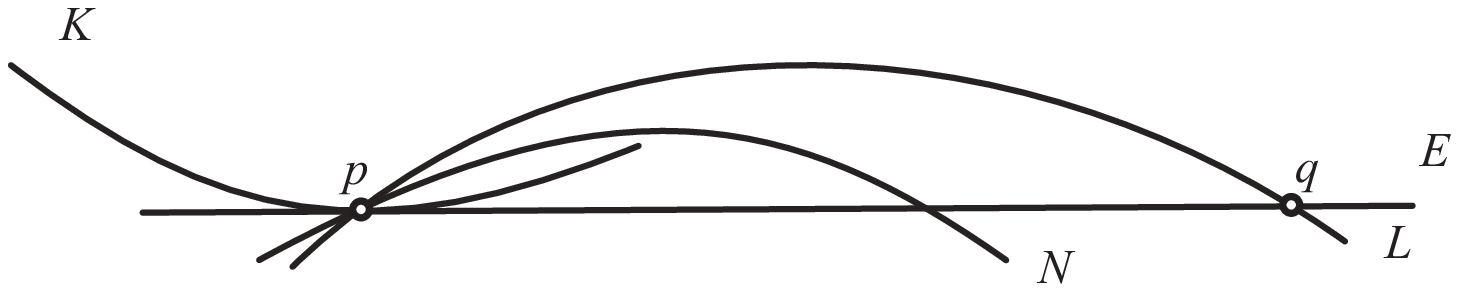}
\caption{}
\end{figure}
\begin{figure}[h]
\includegraphics[width=0.65\textwidth]{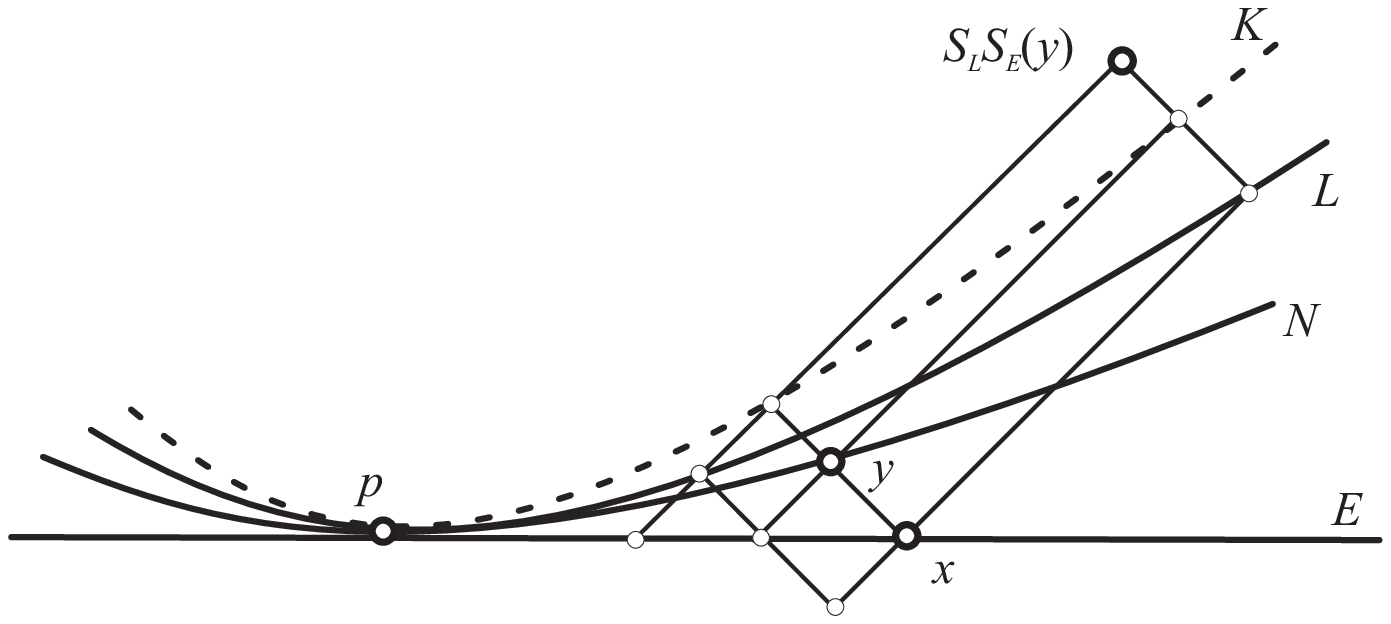}
\caption{}
\end{figure}

\begin{cor}\label{c4.2}
If $E\cap K=\{p\}$, $L\bot E$, $p\in L$, then $i_{L}\circ f_{K}=i_{N}$ for some $N\in\Lambda$ such that $p\in N$.
\end{cor}
\begin{proof}
If $char\mathcal{M}=2$, this is obvious. Let $x\in E$, $y=Kx$, $q=(L(yE))y$ and let $N$ be a circle tangent to $L$ at
$p$ and passing through $q$ (fig. 10). By Proposition \ref{pr4.1}, there exists $K'\in\Lambda$ such that $i_{L}\circ
i_{N}=f_{K'}$ and $K'\cap E=\{p\}$. From the construction of $f_{K'}$ we get $y\in K'$, so $K=K'$.
\end{proof}
\begin{figure}[h]
\includegraphics[width=0.5\textwidth]{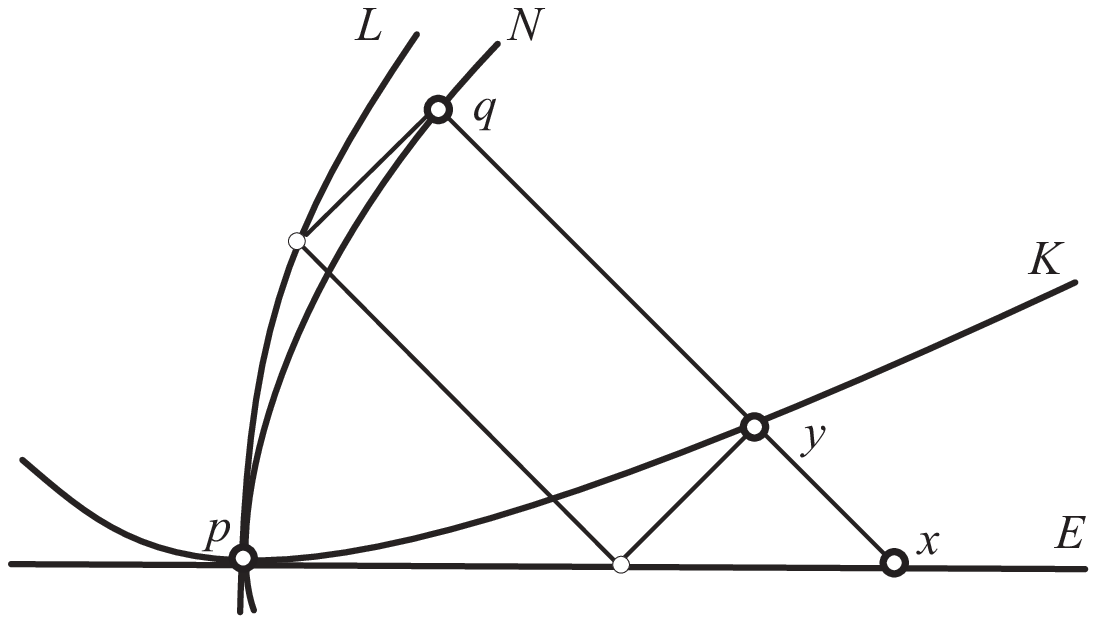}
\caption{}
\end{figure}
\begin{cor}\label{c4.3}
If $p\in E\cap K \cap L\cap M$ and $K,L,M\bot E$, then $i_{M}\circ i_{L}\circ i_{K}=i_{N}$ for some $N\in\Lambda$
through $p$.
\end{cor}
Now we will collect some propositions concerning superpositions of the bijections $f_{K}$. We have not been able to
find a shorter way to prove the postulate (G). The main problem is the case of nonintersecting circles.
\begin{prop}\label{pr4.2}
If $L\bot E$ and $p\in L\cap E$, then for every circle $K$ there exists a circle $M$ such that $f_{K}\circ
i_{L}=f_{M}$.
\end{prop}
\begin{proof}
Assume that $p\notin K$. From Corollary \ref{c4.1}(ii), $f_{K}=i_{S}\circ i_{R}$, where $R,S\bot E$ and $p\in R$,
$Kp\in S$ (fig. 11). If $char\mathcal{M}=2$, then there exists $T\in\Lambda$ such that $p\in T$ and $i_{R}\circ
i_{L}=i_{T}$ by Proposition \ref{pr4.1}. Hence $f_{K}\circ i_{L}=i_{S}\circ i_{R}\circ i_{L}=i_{S}\circ i_{T}$.
According to Corollary \ref{c4.1}(i), there exists a circle $M$ such that $i_{S}\circ i_{T}=f_{M}$. If
$char\mathcal{M}\neq2$, then by Corollary \ref{c4.1}(ii), there exist $T, U\in\Lambda$ such that $i_{S}=i_{U}\circ
i_{T}$, where $p\in T$. By Corollary \ref{c4.3}, we find a circle $W$ such that $i_{W}=i_{T}\circ i_{R}\circ i_{L}$.
Hence $f_{K}\circ i_{L}=i_{S}\circ i_{R}\circ i_{L}=i_{U}\circ i_{T}\circ i_{R}\circ i_{L}=i_{U}\circ i_{W}$. By
Corollary \ref{c4.1}(i), there exists a circle $M$
such that $f_{K}\circ i_{L}=i_{U}\circ i_{W}=f_{M}$.\\
Let now $p\in K$ and $char\mathcal{M}\neq2$. Then we can find $q\neq p$ such that $q\in L\cap E$. If $q\in K$, the
result follows from Corollary \ref{c4.1a}. If $q\notin K$, it follows from the already proven part of the Proposition.\\
Assume $p\in K$ and $char\mathcal{M}=2$. If $K\cap L= \{p\}$, then $f_{K}$ is an involution and the result follows from
Corollary \ref{c4.2}. In the other case let us consider the circle $K':=S_{E}(K)$. We have $K'\cap L=\{p,q\}$ with
$q\notin E$ and we set: $a=Eq$; $r$ such that $r\in K'\cap E$, $r\neq p$; $m=pr$; $n=rp$ (fig. 12). If $b\neq a,p,r$ is
an arbitrary point of $E$, then $f_{K'}=i_{G}\circ i_{F}$, where $F=(m,n,b)^{\circ}$, $G=(m,n,K'b)^{\circ}$ by
Corollary \ref{c4.1}(ii). If $Q$ is a circle such that $f_{Q}=i_{G}\circ i_{L}$ (Corollary \ref{c4.1}), then $b\notin
Q$. Indeed, if $b\in Q$, then $S_{G}(b)\in Q\cap E$ (because $G,L\bot Q,E)$, hence
$E=(S_{G}(b),S_{L}(b),b)^{\circ}=Q$ and we get a contradiction because $n\in Q$ and $n\notin E$. It follows that $b\notin S_{E}(Q)=Q'$.\\
So we have $i_{L}\circ f_{K'}=i_{L}\circ i_{G}\circ i_{F}$, hence $i_{L}\circ f_{K'}=f_{Q}^{-1}\circ i_{F}=f_{Q'}\circ
i_{F}= f_{M'}$ for some circle $M'$ by the first part of the proof. For the circle $M:=S_{E}(M')$, we obtain
$f_{K}\circ i_{L}=f_{M}$.
\end{proof}
\begin{figure}[h]
\includegraphics[width=0.35\textwidth]{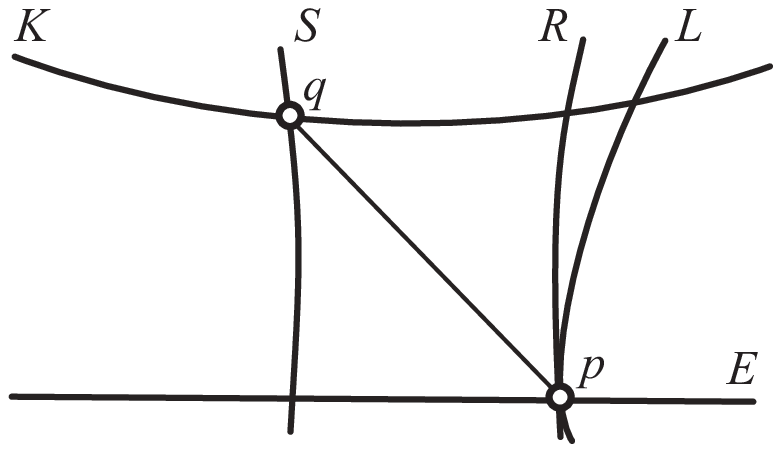}
\caption{}
\end{figure}
\begin{figure}[h]
\includegraphics[width=0.55\textwidth]{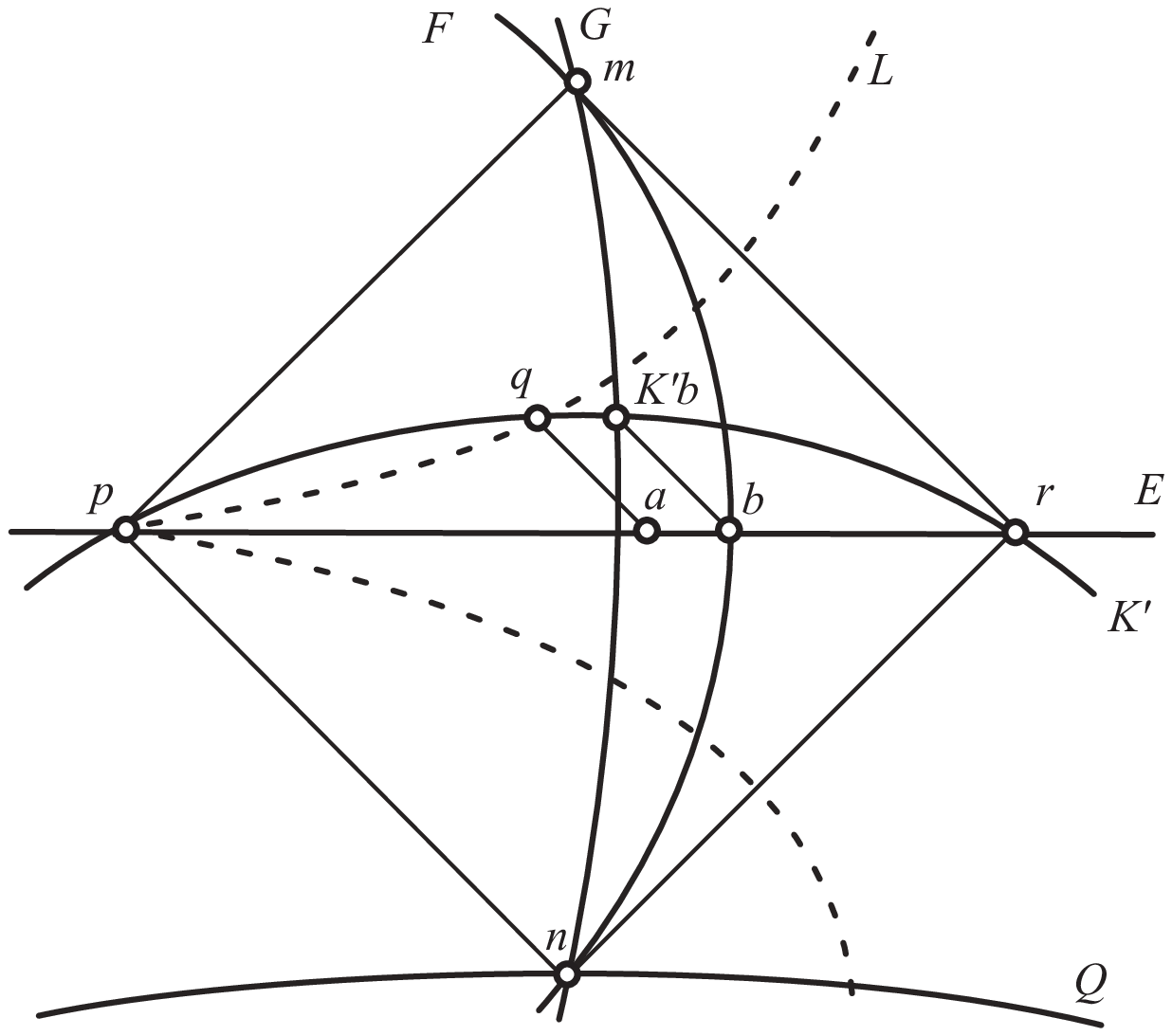}
\caption{}
\end{figure}
\begin{prop}\label{pr4.3}
If $L\bot E$ and $p\in L\cap K$, then there exists a circle $M$ such that $f_{M}=i_{L}\circ f_{K}$.
\end{prop}
\begin{proof}
By Corollary \ref{c4.1}(ii), there exist $R,S\bot E,K$ such that $f_{K}=i_{S}\circ i_{R}$, $p\in S$ and $Ep\in R$. Then
$i_{L}\circ f_{K}=i_{L}\circ i_{S}\circ i_{R}$. As $i_{L}\circ i_{S}=f_{T}$ for some $T$ by Corollary \ref{c4.1a}, we
conclude that $i_{L}\circ f_{K}=f_{T}\circ i_{R}=f_{M}$  for some $M\in\Lambda$ by Proposition \ref{pr4.2}.
\end{proof}
\begin{prop}\label{pr4.4}
If $L,T\bot E$, then there exists $M$ such that $i_{L}\circ i_{T}=f_{M}$.
\end{prop}
\begin{proof}
Let $K$ be an arbitrary circle such that $K\cap L\neq\emptyset$, $K\cap T\neq\emptyset$ and $K\bot T$. Then
$f_{K}=i_{T}\circ i_{S}$ for some $S\in\Lambda$ such that $S\cap E\neq\emptyset$, so we get $i_{L}\circ
i_{T}=(i_{L}\circ f_{K})\circ i_{S}$. By Proposition \ref{pr4.3}, $i_{L}\circ f_{K}=f_{M'}$ for some $M'\in\Lambda$,
hence $i_{L}\circ i_{T}=f_{M'}\circ i_{S}=f_{M}$ for some $M\in\Lambda$ by Proposition \ref{pr4.2}.
\end{proof}
\begin{prop}\label{pr4.5}
If $L\bot E$, then for any $K$ there exists $M$ such that $i_{L}\circ f_{K}=f_{M}$.
\end{prop}
\begin{proof}
By Corollary \ref{c4.1}(ii), there exist circles $P,Q\bot E$ such that $f_{K}=i_{Q}\circ i_{P}$, where $P\cap
E\neq\emptyset$. Hence, $i_{L}\circ f_{K}=(i_{L}\circ i_{Q})\circ i_{P}$ and by Proposition \ref{pr4.2} we can find
$M'$ such that $i_{L}\circ i_{Q}=i_{M'}$. So we have $i_{L}\circ f_{K}=i_{M'}\circ i_{P} = f_{M}$ for some
$M\in\Lambda$ by Proposition \ref{pr4.2}.
\end{proof}
We finally get
\begin{lem}\label{l4.3}
For arbitrary circles $K,L$ there exists a circle $M$ such that $f_{M}=f_{L}\circ f_{K}$.
\end{lem}
Lemma \ref{l4.3} yields
\begin{thm}\label{t4.1}
Every symmetric Minkowski plane satisfies the rectangle axiom \rm(G).
\end{thm}
\begin{cor}\label{c4.4}
A symmetric Minkowski plane satisfies the conclusion of Lemma \ref{l4.1} with assumption:\\
(i') there exists a point $x\in M$ such that $xN\neq xL\wedge K(xN)=(Lx)K$\\
instead (i).
\end{cor}

\end{document}